%  TEXSTUFF +++++++++++++++++++++++++++++++++++++++++++++++++++++++++++++
\input amssym.def
\input epsf

\let \blskip = \baselineskip
\parskip=1.2ex plus .2ex minus .1ex

\tabskip 20pt
\tolerance = 1000
\pretolerance = 50
\newcount\itemnum
\itemnum = 0
\overfullrule = 0pt

% HEADER DEFNS ++++++++++++++++++++++++++++++++++++++++++++++++++++
\def\title#1{\bigskip\centerline{\bigbigbf#1}}
\def\author#1{\bigskip\centerline{\bf #1}\smallskip}
\def\address#1{\centerline{\it#1}}
\def\abstract#1{\vskip1truecm{\narrower\noindent{\bf Abstract.} #1\bigskip}}

%  GENERAL MACROS    +++++++++++++++++++++++++++++++++++++++++++++++++++++
\def\sp{\bigskip}
\def\nosp{\vskip -\the\blskip plus 1pt minus 1pt}

\def\br{\hfil\break} 
\def\ti{\br \hglue \the \parindent}

\def\ce#1{\LP\centerline{#1}}

\def\skipit#1{}
\def\mdag{\raise 3pt\hbox{\dag}}

\def\XP{\par\noindent\hang}
\def\LP{\par\noindent}
\def\BP[#1]{\par\item{[#1]}}
\def\SH#1{\sp\vskip\parskip\leftline{\bigbf #1}\nobreak}

\def\TH#1{\sp\XP{\bf THEOREM\ \shead#1}}
\def\LM#1{\sp\XP{\bf LEMMA\ \shead#1}}

\def\CO#1{\sp\XP{\bf COROLLARY\ \shead#1}}

\def\EX#1{\sp\LP{\bf Example\ \shead#1}}
\def\PF{\LP{\bf Proof:\ }}
\def\NX{\advance\itemnum by 1 \sp\LP {\bf \shead \the\itemnum.\ }}
\def\qed{\null\nobreak\hfill\hbox{${\vrule width 5pt height 6pt}$}\par\sp}

\def\cart{\>\hbox{${\vcenter{\vbox{
    \hrule height 0.4pt\hbox{\vrule width 0.4pt height 4.5pt
    \kern4pt\vrule width 0.4pt}\hrule height 0.4pt}}}$}\>}
\def\bxmu{\>\hbox{${\vcenter{\vbox {
    \hrule height 0.4pt\hbox{\vrule width 0.4pt height 4pt
    \hskip -1.3pt\lower 1.8pt\hbox{$\times$}\negthinspace\vrule width 0.4pt}
    \hrule height 0.4pt}}}$}\>}

\def\lin#1{\hbox to #1true in{\hrulefill}}

% STYLE COMMANDS +++++++++++++++++++++++++++++++++++++++++++++++++++++++

\def\foot#1{\raise 6pt \hbox{#1} \kern -3pt}
% Printed \today.}

% FIGURES AND TABLES +++++++++++++++++++++++++++++++++++++++++++++++++++++++++++
\def\fig #1 #2 #3 #4 #5 {\sp \ce{ {\epsfbox[#1 #2 #3 #4]{figs/#5.ps}} }}
%left x, bottom y, right x, top y, in portrait orientation

%\def\gpic#1{#1 \sp\ce{\box\graph} \medskip} %gpic picture, centered with space
%for troff tables: \s\ccol = centered col, \s\vcol = rule, # to end preamble

% JOURNAL ABBREVS  +++++++++++++++++++++++++++++++++++++++++++++++++++++

\def\JGT{{\it J.\ Graph Theory}}

\def\DAM{{\it Discrete Appl.\ Math.{}}}

\def\ADM{{\it Annals Discr.\ Math.{}}}
\def\GnC{{\it Graphs and Combin.}}

\def\SIDM{{\it SIAM\ J.\ Discr.\ Math.{}}}

\def\JLMS{{\it J.\ Lond. Math.\ Soc.{}}}

\def\CJM{{\it Canad.\ J.\ Math.{}}}
\def\IJM{{\it Israel J.\ Math.{}}}

\def\ANYAS{{\it Annals N.Y.\ Acad.\ Sci.{}}}

% SPECIAL CHARACTERS ++++++++++++++++++++++++++++++++++++++++++++++++++++
\def\al{\alpha}			
\def\eps{\epsilon}

%\def\square{\hbox{${\vcenter{\vbox{
%    \hrule height 0.3pt\hbox{\vrule width 0.3pt height 4.5pt
%    \kern 3.8pt\vrule width 0.3pt}\hrule height 0.3pt}}}$}}

  %{roman {I back 20 N}}
  %{roman {C back 50 C}}
  %{roman {C back 50 Q}}
  %{\rm \kern -8pt Z} }  %\bold Z}   %{roman {Z back 50 Z}}
  %{roman {I back 20 R}}

% BINARY RELATIONS ++++++++++++++++++++++++++++++++++++++++++++++++++++++

		%also \vee,\wedge

    %Note also uppercase doubles line, plus [long][left][right]arrow
    %and [up][down]arrow, plus ne,sw,se,nw+arrow,
    %plus [long]mapsto, hook[left][right]arrow, and harpoons

% DELIMITERS ++++++++++++++++++++++++++++++++++++++++++++++++++++++++++++++
\def\({\left(}	\def\){\right)}

% MACROS WITH ARGUMENTS +++++++++++++++++++++++++++++++++++++++++++++++

\def\CH#1#2{{{#1}\choose{#2}}}

\def\FR#1#2{{#1 \over #2}}

\def\FL#1{\left\lfloor{#1}\right\rfloor}
\def\CL#1{\left\lceil{#1}\right\rceil}

\def\SE#1#2#3{\sum_{#1 = #2} ^ {#3}}

\def\UE#1#2#3{\bigcup_{#1 = #2} ^ {#3}}

\def\VEC#1#2#3{#1_{#2},\ldots,#1_{#3}}

 %"such that"

\def\SET#1:#2{\{#1\colon\;#2\}}

% MODES ++++++++++++++++++++++++++++++++++++++++++++++++++++++++++++++++
\def\B#1{{\bf #1}}		
    %cardinality

% WORDS ++++++++++++++++++++++++++++++++++++++++++++++++++++

% GRAPHS +++++++++++++++++++++++++++++++++++++++++++++++++++++++++++

% POSETS ++++++++++++++++++++++++++++++++++++++++++++++++++++++++++

% MATROIDS +++++++++++++++++++++++++++++++++++++++++++++++++++++++++++

% ENUMERATION ++++++++++++++++++++++++++++++++++++++++++++++++++++++

% PAGE SETUP +++++++++++++++++++++++++++++++++++++++++++++++++++++++++
\magnification=\magstep1
\vsize=9.0 true in
\hsize=6.5 true in
%%\hoffset= -.3 true in
%%\voffset= -.1 true in
\headline={\hfil\ifnum\pageno=1\else\folio\fi\hfil}
%\footline={\hfil-- \folio\ --\hfil}
\footline={\hfil\ifnum\pageno=1\folio\else\fi\hfil}

\parindent=20pt
\baselineskip=12pt
\parskip=.5ex  %1.75

% MACRO STUFF (goes to another file?) ++++++++++++++++++++++++++
\def\shead{ }

\font\bigbf = cmb10 scaled \magstep1

\font\bigbigbf = cmb10 scaled \magstep2

% TITLE PAGE FORMAT ++++++++++++++++++++++++++++++++++++++++++++++++
%\def\titpage{\vbox to \vsize
%     {\vfil \title \vfil \author \vfil \vfil \setabs \vfil \extra}  }

%\input eplain
% gpic macros for use with "gpic -t" as preprocessor for tex,
% by Douglas B. West
 %gpic picture, centered with space

% HEADER SECTION ++++++++++++++++++++++++++++++++++++++++++++++++++++
\title{A NOTE ON GENERALIZED CHROMATIC NUMBER}
\title{AND GENERALIZED GIRTH}
\author{B\'ela Bollob\'as}
\address{University of Memphis and University of Cambridge}
\address{Memphis, TN 38152-6429 and Cambridge CB2 1SB, England}
\author{Douglas B. West\foot{\dag}}
\address{University of Illinois}
\address{Urbana, IL 61801-2975}
\vfootnote{}{\br
   \foot{\dag}Research supported in part by NSA/MSP Grant MDA904-93-H-3040.\br
   Running head: CHROMATIC NUMBER AND GIRTH\br
   Keywords: graph coloring, girth\br
   Written February 1994, revised December 1996 and March 1998.}
\abstract{Erd\H os proved that there are graphs with arbitrarily large
   girth and chromatic number.  We study the extension of this for generalized
   chromatic numbers.}

% DEFINITIONS ++++++++++++++++++++++++++++++++++++++++++++++++++++++++
\def\chp{\chi_H}
\def\chh{\chi_H}

\def\pr{{\rm Prob}}

% DOCUMENT ++++++++++++++++++++++++++++++++++++++++++++++++++++++++++
%\SH
%{1. INTRODUCTION}
Generalized graph coloring describes the partitioning of the vertices into
classes whose induced subgraphs satisfy particular constraints.  When $\B P$ is
a family of graphs, the $\B P$ {\it chromatic number} of a graph $G$, written
$\chi_\B P$, is the minimum size of a partition of $V(G)$ into classes that
induce subgraphs of $G$ belonging to $\B P$.  When $\B P$ is the family of
independent sets, $\chi_\B P$ is the ordinary chromatic number.  General
aspects are studied in [1-3,7-9,11-14,17-18].  Many additional results are known
about particular generalized chromatic numbers.

One aim in the study of generalized chromatic numbers is the extension of
classical coloring results.  Erd\H os\ [4] proved that there exist graphs of
large chromatic number and large girth.  We study the extension of this for
a class of generalized coloring parameters.  We consider the family $\B P$
consisting of all graphs not containing $H$ as a subgraph; we call the
corresponding parameter the $H$-chromatic number and write it as $\chh$.

The natural extension requires an appropriate definition for generalized girth.
For $j\ge2$, an $(H,j)$-{\it cycle} in a graph $G$ is a list of distinct
subgraphs $\VEC H1j$, each isomorphic to $H$, such that $\UE i1j H_i$
contains a cycle that decomposes into $j$ nontrivial paths with the $i$th path
in $H_i$ (any two successive paths in the decomposition share one vertex).  The
$H$-{\it girth} of $G$, written $g_H(G)$, is the minimum $j$ such that $G$
contains an $(H,j)$-cycle, if this exists; otherwise $g_H(G)=\infty$.

One might prefer a weaker notion of cycle.  For $j\ge2$, a {\it weak}
$(H,j)$-{\it cycle} in $G$ is a list of distinct subgraphs $\VEC H1j$, each
isomorphic to $H$, and a selection of distinct vertices $\VEC x1j\in V(G)$ such
that $x_i\in V(H_{i-1})\cap V(H_i)$, with subscripts taken modulo $j$.  The {\it
weak} $H$-girth $g_H^*(G)$ is the minimum $j$ such that $G$ contains a weak
$(H,j)$-cycle, if this exists; otherwise $g_H^*(G)=\infty$.

Every $(H,j)$-cycle is a weak $(H,j)$-cycle, so $g_H^*(G)\le g_H(G)$.
The weak $H$-girth may be considerably smaller than the $H$-girth.
Trivial examples arise when $H$ is disconnected.  Also, when $H$ is
a 3-vertex path and $G=K_{1,3}$, we have $g_H(G)=\infty$ and $g_H^*=2$.

One naturally seeks the existence of graphs with arbitrarily large
$H$-chromatic number and arbitrarily large $H$-girth.  Such a result
does not hold for weak $H$-girth, even when $H$ is $r$-edge-connected.

\EX{}
\it
If $H$ is the union of two copies of $K_{r+1}$ sharing a vertex, then
$\chh(G)\ge 4r+2$ implies $g_H^*(G)\le 2$.
\rm
We prove the contrapositive.  The union of three $r+1$-cliques sharing
a single vertex has weak $H$-girth 2, so $g_H^*(G)>2$ forbids this
as a subgraph.  Thus in $G$ there is no vertex $x$ whose neighborhood induces
a subgraph containing three disjoint $r$-cliques.  Thus in $G[N(x)]$ there
is a set $S_x$ of at most $2r$ vertices (the vertices of a maximal set of
disjoint $r$-cliques) that together contain some vertex of each $r$-clique in
$N(x)$.

Let $G'$ be the spanning subgraph of $G$ whose edges include the edges from
each $x$ to $S_x$ for each $v\in V(G)$.  Since each $S_x$ has size at most
$2r$, each $m$-vertex subgraph of $G'$ has at most $2rm$ edges and thus minimum
degree at most $4r$.  By the Szekeres-Wilf Theorem, $G'$ has ordinary chromatic
number at most $4r+1$.  Every proper ordinary coloring of $G'$ uses at least
two colors on every $r+1$-clique in $G$, since $G'[S]$ has no isolated vertex
when $G[S]$ is an $r+1$-clique.  Thus $\chi_{K_{r+1}}(G)\le 4r+1$.  Since $H$
contains $K_{r+1}$, we have $\chh(G)\le\chi_{K_{r+1}}(G)\le4r+1$.  \qed

Similar examples occur whenever $H$ is not 2-connected, but weak $H$-girth
equals $H$-girth when $H$ is 2-connected.  We prove the desired result using the
stronger concept of $H$-girth.  This also follows from the result of
Erd\H os and Hajnal [5] establishing the existence of $r$-uniform hypergraphs
with large girth and chromatic number (constructive proofs appear in [10,16]).
If $H$ has order $r$, then a graph $G$ with $H$-girth at least $g$ and
$H$-chromatic number $k$ can be obtained from an $r$-uniform hypergraph $\B H$
with girth $g$ and chromatic number $k$ by taking the union of copies of $H$
defined on the vertex sets of the edges of $\B H$.

Erd\H os and Hajnal did not directly define cycles in hypergraphs; instead
they said that an $r$-uniform hypergraph is $s$-{\it circuitless} if for
$t\le s$ every set of $t$ edges contains at least $1+(r-1)t$ vertices in its
union.  When $H$ has $r$ vertices, an $(H,j)$-cycle has at most $(r-1)j$
vertices, so their definition of $s$-circuitless is equivalent to $g_H(G)>s$
when the edges of the hypergraph correspond to copies of $H$ in $G$.

Using the strong version of $H$-girth, we give a short direct existence argument
for graphs with large $H$-chromatic number and large $H$-girth.  When $H$ is a
clique, this becomes a proof of the Erd\H os-Hajnal result.  The main idea of
the construction is similar to theirs, but the computations are somewhat
different, and our presentation is perhaps more self-contained.  When $H$ has
order $r$, the order of our graphs with $g_H>s$ and $\chp>k$ is about $k^{rs}$.
The minimum order of such graphs is addressed in [6], so the point of
our note is its alternative computations.

We summarize the approach.
When seeking $g_H(G)>s$, we say that an $(H,j)$-cycle is a {\it short $H$-cycle}
if $j\le s$.  We will use the ``deletion method'', generating an $n$-vertex 
graph having many copies of $H$ but few short $H$-cycles.  To do this, we let
$r$-sets receive copies of $H$ with probability $p$.  For appropriate $p$, we
expect so many copies of $H$ that in some graph every set of size $\CL{n/k}$
contains more copies of $H$ than the number of short $H$-cycles in the graph.
After deleting edges in copies of $H$ to break all the short $H$-cycles, every
set of size at least $n/k$ still contains a copy of $H$.  By the pigeonhole
principle, the $H$-chromatic number of the resulting graph exceeds $k$.

We need a numerical lemma about tail probabilities in the binomial distribution.

\LM{}
If $X$ has the binomial distribution with $N$ trials and success probability
$p$, then $\pr(X\le pN/2) < 2(2/e)^{pN/2}$.
\PF
For $1\le k\le pN/2$, we have
$${\CH Nk p^k(1-p)^{N-k} \over \CH N{k-1} p^{k-1}(1-p)^{N-k+1}}
= {N-k+1\over k}{p\over 1-p} > 2.$$
In particular, $\pr(X=pN/2-k) <2^{-k}\pr(X=pN/2)$.  Summing, we obtain
$$\pr(X\le pN/2)<2\CH N{pN/2} p^{pN/2}(1-p)^{(1-p/2)N}.$$
Because $\CH N{\al N} < \FR{(1/\al)^{\al N}}{(1-\al)^{(1-\al)N}}$
and $1-\beta<e^{-\beta}$, we conclude
$$\pr(X\le pN/2)<2\cdot2^{pN/2}(\FR{1-p}{1-p/2})^{(1-p/2)N}$$
$$=2\cdot2^{pN/2}(1-\FR{p/2}{1-p/2})^{(1-p/2)N}<2\cdot2^{pN/2}e^{-pN/2}.$$\qed

\TH{}
If $s,k$ are positive
integers with $s>r$, then there is a graph $G$ with $g_H(G)>s$ and $\chp(G)>k$.
Furthermore, if $n$ is sufficiently large, then there is a graph $G$
of order $n$ with $g_H(G)>s$ and $\chp(G)\ge n^{1/(r^\eps s)}$, where $\eps=1$
if $H$ is 2-connected and $\eps=2$ if $H$ is not 2-connected.
\PF
We may assume that $s$ exceeds the order of $H$, since we can add short
cycles to $G$ later to reduce girth.  From vertex set $[n]$ we select
$r$-subsets, independently, each with probability $p$.  This yields a
random $r$-uniform hypergraph $R$.

If $H$ is 2-connected, then set $\ell=s$; otherwise, set $\ell=rs$.
A {\it set-cycle} of length $j$ in $G'$ is a cyclic arrangement of $j$ selected
$r$-sets such that the intersections of successive pairs yield a system of $j$
distinct vertices as representatives.  A set-cycle is {\it short} if it has
length at most $\ell$.  Let $X$ be the number of short $H$-cycles in $R$.
When $n$ is sufficiently large and $N$ and $p$ are appropriately chosen in terms
of $n$, we claim that $E(X)<pN/4$, and that with probability at least $1/2$
every set of $\CL{n/k}$ vertices contains at least $pN/2$ selected $r$-sets.
Thus in some $R$ every $\CL{n/k}$-set contains at least $pN/2$ edges
($r$-sets) of $R$.  Let $R'$ be the hypergraph obtained from $R$ by
deleting some edge ($r$-set) from every short set-cycle in $R$.
On each edge ($r$-set) in $R'$, place a copy of $H$, and let $G$ be the
graph formed by the union of these copies of $H$.

If $H$ is 2-connected, then the only copies of $H$ in $G$ are the ones we
have placed into the $r$-sets of $R'$, so $g_H(G) > \ell=s$.  If $H$ is not
2-connected, then every cycle in $G$ is either entirely in an $r$-set in
$R'$, or else it goes through all the sets of a set-cycle in $R'$.  In the
latter case this cycle has length at least $\ell + 1$, so it is not
contained in the union of $s$ subgraphs isomorphic to $H$.  Hence $g_H(G) > s$
in this case as well.

Finally, the $H$-chromatic number of $G$ is large; since every set of
$\CL{n/k}$ vertices of $R'$ contains an edge of $R'$, we have $\chi_H(G)>k$.

All that remains is to justify the inequalities claimed earlier.  It suffices
to show that our assertions hold when $n$ is sufficiently large and
$k=\FL{n^{1/rs}}$.  Set $N=\CH{n/k}r$ and $p=8(1+\log k)n/(kN)$.
The expected number of $(H,2)$-cycles is less than $\FR{n^2}2 \CH n{r-2}^2 p^2$,
and the expected number of $(H,j)$-cycles is less than
$\FR{n^j}{2j}\CH n{r-2}^j p^j$.  Letting
$\beta= \FR2s \FR{n^{rs-s}p^s}{(r-2)!^s}$, we have
$$\eqalign{E(X) &< \FR{n^2}4 \CH n{r-2} p^2
+\SE j2s \FR{n^j}{2j}\CH n{r-2}^j p^j\cr
&\le \FR{n^s}s \CH n{r-2}^s p^s \le \beta/2\cr}.$$
By Markov's Inequality, the probability now exceeds $1/2$ that $X\le \beta$.

Consider the $\CL{n/k}$-sets.  By the Lemma, the probability that {\it some} 
$\CL{n/k}$-set of vertices contains at most $pN/2$ selected $k$-sets is less
than $\al$, where $\al = 2(2/e)^{pN/2}\CH n{\CL{n/k}}$.  When we have
$k\le n^{1/6}$, we can use Stirling's approximation to conclude that
$\al < \FR12 (2/e)^{pN/2}e^{(1+\log k)n/k}$.
With $pN=8(1+\log k)n/k$, this guarantees $\al \le 1/2$.

Hence there exists an $n$-vertex graph $G$ having at most $pN/2$ selected
$r$-sets and having $X \le \beta$.  It remains only to prove that
$\beta < pN/2 = 4(1+\log k)n/k$.  We need
$p^s\le \FR{2s(1+\log k)}k n^{-s(r-1)+1} (r-2)!^s$.
Using the definition of $p$, it suffices to have
$$16(\log k)k^{r-1}r!n^{-r+1}\le (\FR{2s(1+\log k)}k)^{1/s}(r-2)!n^{-r+1+1/s},$$
or
$$\Bigl(16(\log k)k^{r-1}r(r-1)\Bigr)^s\FR k{2s(1+\log k)} \le n.$$
This is satisfied by the condition on $k$ in terms of $n$.  \qed

Let $H$ be a connected graph.
We say that $G$ is an $H$-{\it forest} if $G$ is a subgraph of a union
of copies of $H$ such that any two of the specified copies of $H$ have
at most one common vertex and every cycle in $G$ is contained in one of
the specified copies of $H$.

\CO{}
Let $H$ be a graph, and let ${\cal F}=\{\VEC F1m\}$ be a family such
that no $F_i$ is an $H$-forest.  Let ${\cal G}$ be the family of graphs 
containing no graph in ${\cal F}$.  For every $k\ge1$, there is a graph
$G\in{\cal G}$ such that $\chi_H(G)=k$.
\PF
If $s$ is the maximum order of the graphs in ${\cal F}$, then every graph
constructed for $s$ in the proof of the theorem belongs to ${\cal G}$.  \qed

In particular, if $m$ is less than the girth of $H$, then $C_m$ is not an
$H$-forest, and hence there exist $C_m$-free graphs with arbitrarily large
$H$-chromatic number.  Also, $\chi_\B P$ is unbounded for $C_m$-free graphs
when $\B P$ is the family of graphs not containing $H$ as an {\it induced}
subgraph.

The Corollary also follows immediately from the result of Ne\v set\v ril and
R\"odl [15] that the class of graphs avoiding a fixed finite set of 2-connected
graphs has the vertex Ramsey property.  They proved their theorem using the
result of Erd\H os and Hajnal [5].

\SH
{REFERENCES}
\BP [1]
M. Borowiecki and P. Mih\'ok, Hereditary properties of graphs,
{\it Advances in Graph Theory} (Vishwa, 1991), 41--68.
\BP [2]
M. Borowiecki, E. Drgas-Burchardt, and P. Mih\'ok, Generalized list colourings
of graphs, {\it Discuss. Math. Graph Theory} 15 (1995), 185--193.
\BP [3]
S.A. Burr and M.S. Jacobson, On inequalities involving vertex-partition
parameters of graphs, {\it Congr.\ Num.{}} 70 (1990), 159--170.
\BP [4]
P. Erd\H os, Graph theory and probability, \CJM\ 11 (1959), 34--38.
\BP [5]
P. Erd\H os and A. Hajnal, On chromatic number of graphs and set-systems,
{\it Acta Math. Acad. Sci. Hung. Tomus} {\bf17}(1966), 61--99.
\BP [6]
P. Erd\H os and L. Lov\'asz, Problems and results on 3-chromatic hypergraphs
and some related questions, {\it Infinite and Finite Sets - Keszthely 1973},
{\it Colloq. Math.  Soc. J\'anos Bolyai} 10 (1975), 
\BP [7]
F. Harary and L.-H.~Hsu, Conditional chromatic number and graph operations,
{\it Bull.\ Inst.\ Math.\ Acad.\ Sinica} 19 (1991), 125--134.
\BP [8]
H. Jacob and H. Meyniel, Extensions of Tur\'an's and Brooks' Theorems and
new notions of stability and coloring in digraphs, \ADM\ 17 (1983), 365--370.
\BP [9]
R.P. Jones, Hereditary properties and P-chromatic numbers, in
{\it Combinatorics} (Proc. Brit. Comb. Conf. 1973), (Cambridge Univ. Press,
1974), 83--88.
\BP [10]
L. Lov\'asz, On chromatic number of graphs and set-systems,
{\it Acta Math. Hungar.} 19 (1968), 59--67.
\BP [11]
D.W. Matula, An extension of Brooks' theorem,
Research Report, Center for Numerical Analysis, Univ. of Texas (Austin, 1973).
\BP [12]
P. Mih\'ok, On vertex partition numbers of graphs, {\it Graphs and other
combinatorial topics (Prague 1982)}, Teubner-Texte Math. 59 (Teubner, 1983),
183--188.
\BP [13]
P. Mih\'ok, On graphs critical with respect to vertex partition numbers.
{\it Discrete Math.} 37 (1981), 123--126.
\BP [14]
P. Mih\'ok, An extension of Brooks' theorem.  {\it 4th Czechoslovakian Symp.
Combinatorics, Graphs, Complexity (Prachatice, 1990)}, {\it Ann. Discrete Math.}
51 (North-Holland, 1992), 235--236.
\BP [15]
J. Ne\v set\v ril and V. R\"odl, Partitions of vertices,
{\it Comm. Math. Univ. Carolinae} 17 (1976), 85--95.
\BP [16]
J. Ne\v set\v ril and V. R\"odl, A short proof of the existence of highly
chromatic hypergraphs without short cycles, {\it J. Combin. Theory Ser. B} 27
(1979), 225--227. 
\BP [17]
E.R. Scheinerman, Generalized chromatic number of random graphs,
\SIDM\ 5 (1992), 74--80.
\BP [18]
M.L. Weaver and D.B. West, Relaxed chromatic numbers of graphs,
\GnC\ 10 (1994), 75--93.

\bye